\documentclass[UTF-8,reqno]{amsart}
\usepackage{enumerate}
\usepackage{amssymb,url,color, booktabs}
\usepackage{mathrsfs}

\numberwithin{equation}{section}

\newcommand{\be}{\begin{eqnarray}}
\newcommand{\ee}{\end{eqnarray}}
\newcommand{\ce}{\begin{eqnarray*}}
\newcommand{\de}{\end{eqnarray*}}
\newtheorem{theorem}{Theorem}[section]
\newtheorem{lemma}[theorem]{Lemma}
\newtheorem{remark}[theorem]{Remark}
\newtheorem{definition}[theorem]{Definition}
\newtheorem{proposition}[theorem]{Proposition}
\newtheorem{Examples}[theorem]{Example}
\newtheorem{corollary}[theorem]{Corollary}
\newtheorem{defi}{Definition}[section]

\usepackage[nobysame]{amsrefs}
\BibSpec{article}{%
+{}{\PrintAuthors} {author}
+{,}{ \textrm} {title}
+{.}{ \textit} {journal}
+{,}{ \textbf} {volume}
+{}{ \parenthesize} {date}
+{,}{ } {pages}
+{.}{ Available at arXiv:} {eprint}
+{.}{} {transition}
}
\BibSpec{book}{%
+{}{\PrintAuthors} {author}
+{,}{ \textit} {title}
+{.}{ \textrm} {series} 
+{,}{ Vol.} {volume}
+{.}{ } {publisher}
+{,}{ } {date}
+{.}{} {transition}
}
\usepackage{color}
\usepackage[colorlinks=true]{hyperref}
\hypersetup{
    linkcolor=blue,          
    citecolor=red,        
    filecolor=blue,      
    urlcolor=cyan
}

\def\e{{\mathrm{e}}}
\def\eps{\varepsilon}

\def\pp{{p}}
\def\p{\partial}

\def\[{{\Big[}}
\def\]{{\Big]}}
\def\<{{\langle}}
\def\>{{\rangle}}
\def\({{\Big(}}
\def\){{\Big)}}

\def\bx{{\mathbf{x}}}

\def\dif{{\mathord{{\rm d}}}}

\def\min{{\mathord{{\rm min}}}}

\def\no{\nonumber}
\def\={&\!\!=\!\!&}
\def\bt{\begin{theorem}}
\def\et{\end{theorem}}
\def\bl{\begin{lemma}}
\def\el{\end{lemma}}
\def\br{\begin{remark}}
\def\er{\end{remark}}

\def\bd{\begin{definition}}
\def\ed{\end{definition}}
\def\bp{\begin{proposition}}
\def\ep{\end{proposition}}
\def\bc{\begin{corollary}}
\def\ec{\end{corollary}}
\def\bx{\begin{Examples}}
\def\ex{\end{Examples}}

\def\cB{{\mathcal B}}

\def\cF{{\mathcal F}}

\def\cM{{\mathcal M}}

\def\cT{{\mathcal T}}

\def\mC{{\mathbb C}}

\def\mE{{\mathbb E}}

\def\mH{{\mathbb H}}
\def\mI{{\mathbb I}}

\def\mL{{\mathbb L}}

\def\mN{{\mathbb N}}

\def\mP{{\mathbb P}}

\def\mR{{\mathbb R}}

\def\1{{\bf 1}}

\def\bP{{\mathbf P}}

\def\sF{{\mathscr F}}

\def\sL{{\mathscr L}}
\def\sM{{\mathscr M}}

\def\sP{{\mathscr P}}

\def\geq{\geqslant}
\def\leq{\leqslant}

\def\div{\mathord{{\rm div}}}

\def\nor{|\mspace{-3mu}|\mspace{-3mu}|}

\def\bP{{\mathbf P}}
\def\bE{{\mathbf E}}

\def\mN{{\mathbb N}}
\def\mR{{\mathbb R}}
\begin{document}

\title{$L^q(L^p)$-theory of stochastic differential equations
}
\author{Pengcheng Xia, Longjie Xie, Xicheng Zhang and Guohuan Zhao}

\address{Pengcheng Xia:
School of Mathematics and Statistics, Wuhan University,
Wuhan, Hubei 430072, P.R.China\\
Email: pcxia@whu.edu.cn
 }

\address{Longjie Xie:
	School of Mathematics and Statistics, Jiangsu Normal University,
	Xuzhou, Jiangsu 221000, P.R.China\\
	Email: longjiexie@jsnu.edu.cn
}

\address{Xicheng Zhang:
School of Mathematics and Statistics, Wuhan University,
Wuhan, Hubei 430072, P.R.China\\
Email: XichengZhang@gmail.com
 }

 \address{Guohuan Zhao:
Fakult\"at f\"ur Mathematik, Universit\"at Bielefeld,
33615, Bielefeld, Germany\\
Email: zhaoguohuan@gmail.com
 }

\thanks{
}

\begin{abstract}
In this paper we show the weak differentiability of the unique strong solution with respect to the starting point $x$
as well as Bismut-Elworthy-Li's derivative formula for the following stochastic differential equation in $\mR^d$:
$$
\dif X_t=b(t,X_t)\dif t+\sigma(t,X_t)\dif W_t,\ \ X_0=x\in\mR^d,
$$
where $\sigma$ is bounded, uniformly continuous and nondegenerate,  $\nabla\sigma\in \widetilde\mL^{p_1}_{q_1}$ and
$b\in \widetilde\mL^{p_2}_{q_2}$ for some $p_i,q_i\in[2,\infty)$ with $\frac{d}{p_i}+\frac{2}{q_i}<1$, $i=1,2$,
where $\widetilde\mL^{p_i}_{q_i}, i=1,2$ are some localized spaces.
Moreover, in the endpoint case $b\in \widetilde\mL^{d; {\rm uni}}_\infty$, we also show the weak well-posedness.

\bigskip
\noindent
\textbf{Keywords}: Krylov's estimate, $L^q(L^p)$-estimates, Zvonkin's transformation, duality.
\\

\thanks{
This work is supported by NNSF grant of China (Nos. 11731009, 11701233), NSF of Jiangsu (No. BK20170226) and the DFG through the CRC 1283
``Taming uncertainty and profiting from randomness and low regularity in analysis, stochastics and their applications''. }

\noindent
{\bf AMS 2010 Mathematics Subject Classification:}  60H10, 60J60.
\end{abstract}

\maketitle \rm

\section{Introduction and main results}

 Consider the following stochastic differential equation (SDE) in $\mR^d$ ($d\geq 2$):
\begin{align}\label{SDE}
\dif X_t=b(t,X_t)\dif t+\sqrt{2}\dif W_t, \ \ X_0=x,
\end{align}
where $(W_t)_{t\geq 0}$ is a $d$-dimensional standard Brownian motion on some filtered probability space $(\Omega,\sF,\bP; (\sF_t)_{t\geq 0})$, and
$b$ is a time-dependent measurable vector field. When $b$ is bounded measurable, Veretennikov \cite{Ve} proved the strong existence and uniqueness of solutions for
SDE \eqref{SDE}. For $T>0$ and $p,q\in(1,\infty)$, let $\mL^p_q(T):=L^q([0,T]; L^p)$. When
$b\in \mL^p_q:=\cap_{T>0}\mL^p_q(T)$ for some $p,q\in[2,\infty)$ with $\frac{d}{p}+\frac{2}{q}<1$, by Girsanov's transformation and
some $\mL^p_q$-estimate for the associated Kolmogorov equation, Krylov and R\"ockner \cite{Kr-Ro} showed the strong well-posedness for SDE \eqref{SDE}
in the class of $X$ that satisfies $\int^T_0|b(t,X_t)|^{2}\dif t<\infty$ a.s. From then on,
there are increasing interests of studying the strong and weak well-posedness for SDE \eqref{SDE} with singular or even distributional drifts, see \cite{Xi-Zh,ZZ} and references therein.

\medskip

After \cite{Kr-Ro}, there are also a lot of works devoted to
studying the properties of the solution $X_t(x,\omega)$ for SDE \eqref{SDE} with singular coefficients.
Among all, we mention that when $b$ is bounded measurable,  Menoukeu etal \cite{M-N-P-Z} showed the weak differentiability of $X_t(x,\omega)$ in $x$
and the Malliavin differentiability of $X_t(x,\omega)$ with respect to the sample point $\omega$.
When $b\in \mL^p_q$  for some $p,q\in[2,\infty)$ with $\frac{d}{p}+\frac{2}{q}<1$ and in the multiplicative noise case, the above regularities in $x$
and $\omega$ were also shown in \cite{Zh4}
by Zvonkin's transformation. However, Zvonkin's transformation used in \cite{Zh4} can not be applied to the bounded drift $b$ because the following PDE
does not allow an $\mH^{2,\infty}$-solution for $b\in L^\infty$ in general:
$$
\p_t u=\Delta u+b\cdot u+b,\ \ u(0)=0.
$$
It should be noticed that the weak differentiability of strong solutions in spatial variables enables us to study the well-posedness of the associated stochastic transport equation
since it is closely related to SDE (\ref{SDE}) through the stochastic inverse flow induced by the strong solution, see \cites{F-G-P,Mo-Ni-Pr}
and references therein.
One of the aim of this paper is to provide a unified treatment for the main results in \cite{M-N-P-Z} and \cite{Zh4}
and extends them to the case of {\it local} integrable coefficients.

\medskip

On the other hand, in the critical case $\frac{d}{p}+\frac{2}{q}=1$ with $p,q\in[2,\infty)$,
Beck etal \cite{Be-Fl-Gu-Ma} claimed the existence and uniqueness of strong solutions to SDE \eqref{SDE} for {\it almost all} starting point $x$.
Recently, when $b$ belongs to some Lorentz space $L^{q,1}(L^p)\subset L^{q,q}(L^p)=\mL^p_q$
for some $p,q\in[2,\infty)$ with $\frac{d}{p}+\frac{2}{q}=1$, still by Zvonkin's transformation,
Nam \cite{Na} showed the existence and uniqueness of strong solutions for SDE \eqref{SDE}.
When $b\in L^d(\mR^d)$ is {\it time-independent}, Kinzebulatov and Semenov \cite{Ki-Se} showed the existence of weak solutions
for each starting point $x\in\mR^d$, but the uniqueness is left open.
Moreover, in the supercritical case $b\in \mL^p_q$ for some $p,q\in[2,\infty)$ with $\frac{d}{p}+\frac{2}{q}<2$,
under an extra integrability assumption on $(\div b)^-$, in a recent work \cite{Zh-Zh2}, the last two authors of the present paper showed the existence of weak solutions.
Another goal of this paper is to show the existence and uniqueness of weak solutions for SDE \eqref{SDE} with multiplicative noise
in the endpoint case $b\in \widetilde\mL^{d; {\rm uni}}_\infty$, which is not covered by all of the above results.

\medskip

In this paper, we shall consider the following SDE driven by multiplicative Brownian noises:
\begin{align}\label{SDE1}
\dif X_t=b(t,X_t)\dif t+\sigma(t,X_t)\dif W_t, \ \ X_0=x,
\end{align}
where $\sigma:\mR_+\times\mR^d\to\mR^d\otimes\mR^d$ and $b:\mR_+\times\mR^d\to\mR^d$ are Borel measurable functions. The generator of this SDE is given by
\begin{align}\label{LR2}
\sL^{\sigma,b}_t f(x):=\tfrac{1}{2}(\sigma^{ik}\sigma^{jk})(t,x)\p_{i}\p_j f(x)+b^i(t,x)\p_if(x).
\end{align}
Here and below, we use Einstein's convention that the repeated indices in a product will be summed automatically. Throughout this paper, we assume that
\begin{enumerate}[{\bf (H$^\sigma$)}]
\item $\lim_{|x-y|\to 0}\sup_t\|\sigma(t,x)-\sigma(t,y)\|_{HS}=0$,
and for some $c_0\geq 1$ and  for all $(t,x)\in\mR_+\times\mR^d$,
\begin{align*}
c_0^{-1}|\xi|^2\leq |\sigma(t,x)\xi|^2\leq c_0|\xi|^2,\ \ \forall\xi\in\mR^d,
\end{align*}
where $\|\cdot\|_{HS}$ stands for the Hilbert-Schmidt norm of a matrix.
\end{enumerate}

Our first main result in this paper is:
\bt\label{Main1}
Assume {\bf (H$^\sigma$)} and $\nabla\sigma\in \widetilde\mL^{p_1}_{q_1}, b\in \widetilde\mL^{p_2}_{q_2}$
for some $p_i,q_i\in[2,\infty)$ with $\frac{d}{p_i}+\frac{2}{q_i}<1$, $i=1,2$,
where $\widetilde\mL^{p}_{q}$ is defined by \eqref{GG1} below.
Then for each $x\in\mR^d$, there is a unique strong solution $X_t(x)$ for SDE \eqref{SDE1}.
Moreover, $X_t(x)$ enjoys the following properties:
\begin{enumerate}[(i)]
\item(Krylov's estimate) For any $p,q\in(1,\infty)$ with $\frac{d}{p}+\frac{2}{q}<2$ and $T>0$,
there is a constant $C>0$ such that for all $x\in\mR^d$ and $0\leq t_0<t_1\leq T$, $f\in\widetilde\mL^p_q(t_0,t_1)$,
\begin{align*}
\bE\left(\int^{t_1}_{t_0}f(s,X_s(x))\dif s\Big|\sF_{t_0}\right)\leq C\nor f\nor_{\widetilde\mL^p_q(t_0,t_1)},
\end{align*}
where $\nor \cdot\nor_{\widetilde\mL^p_q(t_0,t_1)}$ is defined by \eqref{GG1} below.

\item(Weak differentiability) For each $t\geq 0$, the mapping $x\mapsto X_t(x)$ is almost surely weak differentiable and for any $T>0$ and
$p\geq 1$,
\begin{align}\label{Gr}
\sup_{x\in \mR^d}\bE\left(\sup_{t\in[0,T]}|\nabla X_t(x)|^p\right)<\infty.
\end{align}
\item(Derivative formula) For any $t>0$ and $\varphi\in C^1_b(\mR^d)$, it holds that for Lebesgue-almost all $x\in\mR^d$,
\begin{align}
\nabla\bE \varphi(X_t(x))=\frac{1}{t}\bE\left(\varphi(X_t(x))\int^t_0\sigma^{-1}(s,X_s(x))\nabla X_s(x)\dif W_s\right).\label{feller}
\end{align}
\end{enumerate}
\et
\br
As we mentioned before, when $\nabla\sigma,b\in\mL^{p}_{q}$ for some $p,q\in(2,\infty)$ with $\frac{d}{p}+\frac{2}{q}<1$, the above theorem has been obtained in \cite{Zh4}.
Notice that $b\in\mL^\infty$ is not covered by \cite{Zh4}.
The novelty of our result here is that we are considering some localized $\widetilde\mL^{p}_{q}$-spaces so that we still have the global properties \eqref{Gr} and \eqref{feller}.
In particular, we extend the main results in \cite{M-N-P-Z, Mo-Ni-Pr, Zh4} to more general cases, and our proofs are much simpler than \cite{M-N-P-Z}.
\er

Let $\mC$ be the space of all continuous functions from $\mR_+$ to $\mR^d$ endowed with the usual Borel $\sigma$-field $\cB(\mC)$, and $\omega_t$ 
the canonical process over $\mC$. For $t\geq 0$, let $\cB_t:=\cB_t(\mC)$ be the natural filtration generated by $\{\omega_s: s\leq t\}$.
All the probability measures over $(\mC,\cB(\mC))$ is denoted by $\sP(\mC)$.
We introduce the following notion of martingale solutions.
\bd
Given $(s,x)\in\mR_+\times\mR^d$, we call a probability measure $\mP_{s,x}\in\sP(\mC)$ a martingale solution of SDE \eqref{SDE1} with starting point $(s,x)$ if
$\mP_{s,x}(\omega_t=x, t\leq s)=1$, and for all $f\in C^2_b(\mR^d)$, $M^f_t$ is a $\cB_t$-martingale under $\mP_{s,x}$, where
$$
M^f_t(\omega):=f(\omega_t)-f(x)-\int^t_s\sL^{\sigma,b}_r f(\omega_r)\dif r,\ \ t\geq s,
$$
and $\sL^{\sigma,b}_r$ is defined by \eqref{LR2}.
All the martingale solution $\mP_{s,x}$ of SDE \eqref{SDE1} with starting point $(s,x)$ and coefficients $(\sigma,b)$ is denoted by $\sM^{\sigma,b}_{s,x}$.
\ed

Our second main result is the following weak well-posedness of SDE (\ref{SDE1}) in the endpoint case $b\in \widetilde\mL^{d; {\rm uni}}_\infty$
(see \eqref{GR1} below for the definition of $\widetilde\mL^{d; {\rm uni}}_\infty$).

\bt\label{Main2}
Assume {\bf (H$^\sigma$)} holds and $b\in \widetilde\mL^{d; {\rm uni}}_\infty$.
Then for each $(s,x)\in\mR_+\times\mR^d$, there is a
 unique martingale solution $\mP_{s,x}\in\sM^{\sigma,b}_{s,x}$ for SDE \eqref{SDE1} which satisfies
that for any $p,q\in(1,\infty)$ with $\frac{d}{p}+\frac{2}{q}<2$ and $T>0$,
there is a constant $C>0$ such that for all $x\in\mR^d$ and $s\leq t_0<t_1\leq T$,  $f\in\widetilde\mL^p_q(t_0,t_1)$,
\begin{align}\label{KR1}
\mE^{\mP_{s,x}}\left(\int^{t_1}_{t_0}f(r,\omega_r)\dif r\Big|\cB_{t_0}\right)\leq C\nor f\nor_{\widetilde\mL^p_q(t_0,t_1)}.
\end{align}
\et

The proof of our main results relies on the
$\mL^p_q$-maximal regularity estimate for the following second order parabolic PDE in $\mR_+\times\mR^d$:
\begin{align}\label{PDE101}
\partial_{t}u=a^{i j}\p_{i}\p_ju+f,\ \ u(0)=0,
\end{align}
where $a(t,x):\mR_+\times\mR^d\to\mR^d\otimes\mR^d$ is a symmetric matrix-valued  Borel function and satisfies
\begin{enumerate}[{\bf (H$^a$)}]
\item $\lim_{|x-y|\to 0}\sup_{t\in\mR_+}\|a(t,x)-a(t,y)\|_{HS}=0$ and for some $c_0\geq 1$
and for all $(t,x)\in\mR_+\times\mR^d$,
\begin{align}\label{Non}
c_0^{-1}|\xi|^2\leq a^{ij}(t,x)\xi_i\xi_j\leq c_0|\xi|^2,\ \ \forall\xi\in\mR^d.
\end{align}
\end{enumerate}
More precisely, for any $p,q\in(1,\infty)$, we want to establish the following estimate:
\begin{align}\label{KD1}
\|\p_t u\|_{\mL^p_q(T)}+\|\nabla^2 u\|_{\mL^p_q(T)}\leq C\|f\|_{\mL^p_q(T)}.
\end{align}
Such type of estimate has been used in \cite{Xi-Zh} to study
 the strong well-posedness of SDEs with Sobolev diffusion coefficients.
Notice that when $p=q$, it is a standard procedure to prove \eqref{KD1} by freezing coefficient argument (cf. \cite{Zh4}).
While for $p\not=q$, it is non-trivial. When $a^{ij}$ is independent of $x$, \eqref{KD1} was first proved by Krylov in \cite{K3}.
In the spatial dependent case, Kim \cite{Ki} showed
\eqref{KD1} only for $p\leq q$. Here we shall drop this restriction by a duality method.
In particular, we need to treat the adjoint equation of \eqref{PDE101} in Sobolev spaces with negative differentiability index, 
see Theorem \ref{Th32} below, which is of independent interest. Moreover, we also show the estimate \eqref{KD1} in localized space $\widetilde\mL^p_q(T)$.

\medskip

This paper is organized as follows: In Section 2, we collect some preliminary tools.
Section 3 is devoted to the study of $\mL^p_q$-maximal regularity estimate for second order parabolic equations. In Section 4, we prove our main theorems.
Throughout this paper we shall use the following conventions:
\begin{itemize}
\item The letter $C$ denotes a constant, whose value may change in different places.
\item We use $A\lesssim B$ and $A\asymp B$  to denote $A\leq C B$ and $C^{-1} B\leq A\leq CB$ for some unimportant constant $C>0$, respectively.
\item For any $\eps\in(0,1)$, we use $A\lesssim\eps B+D$ to denote
$A\leq \eps B+C_\eps D$ for some constant $C_\eps>0$.
\item $\mN_0:=\mN\cup\{0\}$, $\mR_+:=[0,\infty)$, $a\vee b:=\max(a,b)$, $a\wedge b:=\min(a,b)$, $a^+:=a\vee 0$.
\item $\nabla_x:=\p_x:=(\p_{x_1},\cdots,\p_{x_d})$, $\p_i:=\p_{x_i}:=\p/\p x_i$.
\end{itemize}

\section{Preliminaries}

First of all, we  introduce some spaces and notations for later use.
For $(\alpha,p)\in\mR\times(1,\infty)$,  let $H^{\alpha,p}:=(\mI-\Delta)^{-\alpha/2}\big(L^p(\mR^d)\big)$
be the usual Bessel potential space with norm
$$
\|f\|_{\alpha,p}:=\|(\mI-\Delta)^{\alpha/2}f\|_p,
$$
where $\|\cdot\|_p$ is the usual $L^p$-norm in $\mR^d$, and $(\mI-\Delta)^{\alpha/2}f$ is defined through Fourier's  transform
$$
(\mI-\Delta)^{\alpha/2}f:=\cF^{-1}\big((1+|\cdot|^2)^{\alpha/2}\cF f\big).
$$
Notice that for $n\in\mN$ and $p\in(1,\infty)$, an equivalent norm in $H^{n,p}$ is given by
$$
\|f\|_{n,p}=\|f\|_p+\|\nabla^n f\|_{p}.
$$
Let $\chi\in C^\infty_c(\mR^d)$ be a smooth function with $\chi(x)=1$ for $|x|\leq 1$ and $\chi(x)=0$ for $|x|>2$.
For $r>0$ and $z\in\mR^d$, define
\begin{align}\label{CHI}
\chi_r(x):=\chi(x/r),\ \ \chi^z_r(x):=\chi_r(x-z).
\end{align}
Fix $r>0$. We introduce the following localized $H^{\alpha,p}$-space:
$$
\widetilde H^{\alpha,p}:=\Big\{f\in H^{\alpha,p}_{loc}(\mR^d),\nor f\nor_{\alpha,p}:=\sup_z\|\chi^z_r f\|_{\alpha,p}<\infty\Big\}.
$$
For $T>0$, $p,q\in(1,\infty)$ and $\alpha\in\mR$, we also define space-time function space
$$
\mL^p_q(T):=L^q\big([0,T];L^p\big),\ \  \mH^{\alpha,p}_q(T):=L^q\big([0,T];H^{\alpha,p}\big),
$$
and  the localized space $\widetilde\mH^{\alpha,p}_q(T)$ with norm
\begin{align}\label{GG1}
\nor f\nor_{\widetilde\mH^{\alpha,p}_q(T)}:=\sup_{z\in\mR^d}\|\chi^z_r f\|_{\mH^{\alpha,p}_q(T)}<\infty.
\end{align}
For $q=\infty$ and $p\in[1,\infty)$, we define $\widetilde\mL^{p; {\rm uni}}_\infty(T)$ being all the functions $f\in\widetilde\mL^{p}_\infty(T)$ with
\begin{align}\label{GR1}
\lim_{\eps\to 0}\sup_{t\in[0,T]}\nor f(t,\cdot)*\rho_\eps-f(t,\cdot)\nor_p=:\lim_{\eps\to 0}\kappa^f_T(\eps)=0,
\end{align}
where $(\rho_\eps)_{\eps\in(0,1)}$ is a family of mollifiers in $\mR^d$. 
For simplicity we shall write
$$
H^{\infty,p}:=\cap_{\alpha>0} H^{\alpha,p},\ \ 
\widetilde\mH^{\alpha,p}_q:=\cap_{T>0}\widetilde\mH^{\alpha,p}_q(T),\ \ \widetilde\mL^{p}_q:=\cap_{T>0}\widetilde\mL^{p}_q(T).
$$
It is not hard to show that the definitions of $\widetilde H^{\alpha,p}$ and $\widetilde\mH^{\alpha,p}_q(T)$ do not depend on the choice of $r$ and $\chi$.
In fact, we can prove that for any $r,r'>0$ (cf. \cite{Zh-Zh2}),
\begin{align}\label{GW2}
\sup_{z\in\mR^d}\|\chi^z_r f\|_{\mH^{\alpha,p}_q(T)}\asymp\sup_{z\in\mR^d}\|\chi^z_{r'} f\|_{\mH^{\alpha,p}_q(T)}.
\end{align}
Notice that
$$
L^q([0,T];\widetilde H^{\alpha,p})\subset\widetilde\mH^{\alpha,p}_q(T).
$$
Now we list some easy properties about space $\widetilde\mH^{\alpha,p}_q(T)$ for later use.
\begin{itemize}
\item The following Sobolev embedding holds: For any $\alpha>0$ , $p,q\in[1,\infty)$ and $p'\in[p,\tfrac{pd}{d-p\alpha}\1_{p\alpha<d}+\infty\cdot\1_{p\alpha>d}]$, there is a constant $C>0$ such that
\begin{align}\label{Sob}
\nor f\nor_{\widetilde\mL^{p'}_q(T)}\leq C\nor f\nor_{\widetilde\mH^{\alpha,p}_q(T)}.
\end{align}

\item For any $f\in\widetilde\mH^{\alpha,p}_q$, it holds that for any $T, R>0$ (cf. \cite[Proposition 4.1]{Zh-Zh2}),
\begin{align}\label{LQ1}
\sup_\eps\nor f_\eps\nor_{\widetilde\mH^{\alpha,p}_q(T)}\leq C\nor f\nor_{\widetilde\mH^{\alpha,p}_q(T)},\ \lim_{\eps\to 0}\nor (f_\eps-f)\chi_R\nor_{\widetilde\mH^{\alpha,p}_q(T)}=0,
\end{align}
where $f_\eps:=f*\rho_\eps$ is the usual mollifying approximation of $f$.

\item Let $p,q\in[2,\infty)$ satisfy $\frac{d}{p}+\frac{2}{q}<2$.
If $u\in\widetilde\mH^{2,p}_q(T)$ and $\p_tu\in\widetilde\mL^{p}_q(T)$, then $u\in C([0,T]\times\mR^d)$
(cf. \cite[Lemma 10.2]{Kr-Ro}).
\end{itemize}

For $R\in(0,\infty)$, we define the local Hardy-Littlewood maximal function by
$$
\cM_Rf(x):=\sup_{r\in(0,R)}\frac{1}{|B_r|}\int_{B_r}f(x+y)\dif y,
$$
where $B_r:=\{x\in\mR^d: |x|<r\}$ is the ball in $\mR^d$.
We have the following  results (cf. \cite{St} or \cite[]{Zh1}).
\bl\label{Le2}
(i) For any $R>0$, there exists a constant $C=C(d,R)>0$ such that for any $f\in L^\infty(\mR^d)$ with $\nabla f\in L^1_{loc}(\mR^d)$
and Lebesgue-almost all $x,y\in \mR^d$,
\begin{align}
|f(x)-f(y)|\leq C |x-y|(\cM_R|\nabla f|(x)+\cM_R|\nabla f|(y)+\|f\|_\infty).\label{ES2}
\end{align}
(ii) For any $p>1$, $q\geq 1$ and $R>0$, there is a constant $C=C(R,d,p)>0$ such that for all $f\in \widetilde\mL^p_q(T)$,
\begin{align}\label{GW1}
\nor\cM_R f\nor_{\widetilde\mL^p_q(T)}\leq C\nor f\nor_{\widetilde\mL^p_q(T)}.
\end{align}
\el
\begin{proof}
(i) If $|x-y|\leq R$, then by \cite[Lemma 5.4]{Zh1} we have
$$
|f(x)-f(y)|\leq C |x-y|(\cM_R|\nabla f|(x)+\cM_R|\nabla f|(y)).
$$
If $|x-y|>R$, then
$$
|f(x)-f(y)|\leq 2|x-y|\,\|f\|_\infty/R.
$$
Thus (\ref{ES2}) is true.
\medskip\\
(ii) Noticing that for $|y|\leq R$, $\chi_R(x)=\chi_R(x)\chi_{3R}(x+y)$, by definition we have
\begin{align*}
\|\chi^z_{R}\cM_R f_s\|^p_p&=\int_{\mR^d}\left|\chi_{R}(x)\sup_{r\in(0,R)}\frac{1}{|B_r|}\int_{B_r}f_s(x+z+y)\dif y\right|^p\dif x\\
&\leq\int_{\mR^d}\left(\sup_{r\in(0,R)}\frac{1}{|B_r|}\int_{B_r}\chi_{3R}(x+y)|f_s|(x+z+y)|\dif y\right)^p\dif x\\
&\leq C\|\chi_{3R}\cdot f_s(\cdot+z)\|^p_p=C\|\chi^z_{3R} f_s\|^p_p,
\end{align*}
which in turn gives \eqref{GW1} by \eqref{GW2}.
\end{proof}

The following freezing lemma is taken from \cite[Lemma 4.1]{ZZ}.
\bl\label{Le21}
Let $\phi$ be a nonzero smooth function with compact support.
Define $\phi_z(x):=\phi(x-z)$. For any $\alpha\in\mR$ and $p\in (1,\infty)$, there exists a constant $C\geq 1$ depending only on $\alpha,p,\phi$ such that
for all $f\in H^{\alpha,p}$,
\begin{align}\label{GJ1}
C^{-1} \|f\|_{\alpha,p}\leq \left(\int_{\mR^d}  \|\phi_z f\|_{\alpha,p}^p\dif z\right)^{1/p} \leq C \|f\|_{\alpha,p}.
\end{align}
\el

The following lemma was proven in \cite{K3} (see also \cite[Lemma 2.5]{Ki}).
\bl\label{Le01}
For $k=1,\cdots,n$, let $a_k:\mR\to\mR^d\otimes\mR^d$ be a measurable function and satisfy that for some $c_0\geq 1$,
$$
c_0^{-1}|\xi|^2\leq a^{ij}_k(t)\xi_i\xi_j\leq c_0|\xi|^2,\ \ \forall(t,\xi)\in\mR\times\mR^d,
$$
For fixed $\alpha\in\mR$, $p\in(1,\infty)$ and $\lambda\geq 0$, let $u_k\in \mH^{\alpha,p}_p$ solve the following PDE in the distributional sense:
$$
\p_t u_k=a^{ij}_k\p_{ij} u_k-\lambda u_k+f_k,\ \ u(0)=0.
$$
Then for any $T\geq 0$, there is a constant $N=N(d,\alpha,p,n,c_0)>0$ independent of $T,\lambda$ such that
$$
\int^T_0\prod_{k=1}^n\|\nabla^2 u_k(t)\|_{\alpha,p}^p\dif t\leq N\sum_{k=1}^n\int^T_0\|f_k\|_{\alpha,p}^p\prod_{\ell\not=k}\|\nabla^2 u_\ell(t)\|_{\alpha,p}^p\dif t.
$$
\el

\section{$\widetilde\mL^p_q$-maximal regularity estimate for parabolic equations}

Consider the following second order parabolic PDE in $\mR_+\times\mR^d$:
\begin{align}\label{PDE}
\partial_{t}u=a^{i j}\p_i\p_ju+b^i\p_i u-\lambda u+f,\ \ u(0)=0,
\end{align}
where $\lambda\geq 0$, $a(t,x):\mR_+\times\mR^d\to\mR^d\otimes\mR^d$ and $b(t,x):\mR_+\times\mR^d\to\mR^d$ are Borel measurable functions.
The main aim of this section is to establish the following $\widetilde\mL^p_q$-maximal regularity estimate for the above equation.

\bt\label{Th1}
Let $p,q\in(1,\infty)$. Assume {\bf (H$^a$)} and one of the following conditions holds:
\begin{enumerate}[(i)]
\item{(Subcritical case)} $\frac{d}{p}+\frac{2}{q}<1$ and for any $T>0$, $\nor b\nor_{\widetilde\mL^{p}_{q}(T)}\leq\kappa^b_T<\infty$;
\item{(Critical case)} $p\in(1,d)$ and $b\in \widetilde\mL^{d; {\rm uni}}_\infty$.
\end{enumerate}
Then for any $f\in\widetilde\mL^p_q$ and $\lambda\geq 1$, there exists a unique strong solution $u\in\widetilde \mH^{2,p}_q$
to PDE \eqref{PDE}, that is, for all $t\geq 0$ and Lebesgue almost all $x\in\mR^d$,
$$
u(t,x)=\int^t_0 (a^{ij}\p_i\p_j)u(s,x)\dif s+\int^t_0 (b^i\p_iu)(s,x)\dif s-\lambda\int^t_0u(s,x)\dif s+\int^t_0 f(s,x)\dif s.
$$
Moreover, for any $T>0$ and $\alpha\in[0,2-\frac{2}{q})$,
there is a constant $C>0$ only depending on $\alpha,p,q,d,c_0,T$ and the continuity modulus of $a$, as well as $\kappa^b_T$ in case (i),
and $\kappa^b_T(\eps)$ in case (ii), where $\kappa^b_T(\eps)$ is defined by \eqref{GR1}, such that for any  $\lambda\geq1$,
\begin{align}\label{Max}
\lambda^{1-\frac{\alpha}{2}-\frac{1}{q}}\nor u\nor_{\widetilde\mH^{\alpha,p}_\infty(T)}
+\nor \p_t u\nor_{\widetilde\mL^{p}_q(T)}+\nor u\nor_{\widetilde\mH^{2,p}_q(T)}\leq C\nor f\nor_{\widetilde\mL^p_q(T)}.
\end{align}
\et

\br
In critical case (ii), if $b(t,x)=b(x)\in L^d(\mR^d)$ is time-independent, then $b\in  \widetilde\mL^{d; {\rm uni}}_\infty$.
\er

\subsection{Smooth $a$ and $f$}
In this subsection, we study PDE (\ref{PDE}) with $b\equiv 0$ and  $a$ smooth enough, that is,  $a$ satisfies {\bf (H$^a$)} and for all $m\in\mN$,
$$
\|\nabla^m a^{ij}\|_\infty<\infty,
$$
where $\nabla^m$ stands for the $m$-order gradient.
Given $s<t$, $\lambda\geq 0$ and $\varphi,\psi\in C^\infty_b(\mR^d)$,  consider the following forward heat equation
\begin{align}\label{Eq1}
\p_t u=a^{ij}\p_{ij} u-\lambda u,\ u(s)=\varphi,
\end{align}
and backward (adjoint) heat equation
\begin{align}\label{Eq2}
\p_s w=\lambda w-\p_{ij} (a^{ij}w),\ w(t)=\psi.
\end{align}
Let $u(t)$ and $w(s)$ be the unique solutions of \eqref{Eq1} and \eqref{Eq2} respectively. We shall simply write
$$
\cT_{s,t}\varphi:=u(t),\ \ \cT^*_{s,t}\psi:=w(s).
$$
In other words, we have
$$
\p_t\cT_{s,t}\varphi=a^{ij}\p_{ij}\cT_{s,t}\varphi-\lambda \cT_{s,t}\varphi,\ \
\p_s\cT^*_{s,t}\psi=\lambda \cT^*_{s,t}\psi-\p_{ij}(a^{ij}\cT^*_{s,t}\psi).
$$
Le $p\geq 1$. By the chain rule and above equations, it is easy to see that for any $\varphi,\psi\in H^{\infty,p}\subset C^\infty_b(\mR^d)$,
$$
\<\cT_{s,t}\varphi,\psi\>-\<\varphi,\cT^*_{s,t}\psi\>=\int^t_s\dif_r\<\cT_{s,r}\varphi,\cT^*_{r,t}\psi\>=0,
$$
where $\<f,g\>:=\int_{\mR^d}f(x)g(x)\dif x$, which means that
\begin{align}\label{Dual}
\<\cT_{s,t}\varphi,\psi\>=\<\varphi,\cT^*_{s,t}\psi\>.
\end{align}
Fix $T>0$ and $p,q\geq 1$. For $f\in \mL^q_T(H^{\infty,p}):=L^q([0,T]; H^{\infty,p})$, define
\begin{align}\label{Eq6}
u(t,x):=\int^t_0\cT_{s,t}f(s,x)\dif s,\ \ w(s,x):=\int^T_s\cT^*_{s,t}f(t,x)\dif t.
\end{align}
It is well known that $u$ solves the following forward equation
\begin{align}\label{PDE11}
\partial_{t}u=a^{i j}\p_{ij}u-\lambda u+f,\ \ u(t)|_{t\leq 0}=0,
\end{align}
and $w$ solves the following backward equation
\begin{align}\label{PDE12}
\partial_sw=\lambda w-\p_{ij}(a^{i j}w)-f,\ \ w(s)|_{s\geq T}=0.
\end{align}

We first prove the following a priori estimates by duality.

\bt\label{Th32}
Under {\bf (H$^a$)},
for any $p,q\in(1,\infty)$ and $T>0$, there is a constant $C>0$ only depending on $T, d,p,q,c_0$ and the continuity modulus of $a$ such that
for any $f\in \mL^q_T(H^{\infty,p})$ and $\lambda\geq 0$,
\begin{align}
\|\nabla^2 u_\lambda\|_{\mL^p_q(T)}&\leq C\|f\|_{\mL^p_q(T)},\label{Eq3}\\
\|\nabla^2 w_\lambda\|_{\mH^{-2,p}_q(T)}&\leq C \|f\|_{\mH^{-2,p}_q(T)},\label{Eq39}
\end{align}
where $u_\lambda$ and $w_\lambda$ are solutions of \eqref{PDE11} and \eqref{PDE12}, respectively.
Moreover, for any $\alpha\in[0,2-\frac{2}{q})$, we also have
\begin{align}
\|u_\lambda\|_{\mH^{\alpha,p}_\infty(T)}&\leq C(1\vee\lambda)^{\frac{\alpha}{2}-1+\frac{1}{q}}\|f\|_{\mL^p_q(T)},\label{Eq30}\\
\|w_\lambda\|_{\mH^{\alpha-2,p}_\infty(T)}&\leq C (1\vee\lambda)^{\frac{\alpha}{2}-1+\frac{1}{q}}\|f\|_{\mH^{-2,p}_q(T)}.\label{Eq301}
\end{align}
\et
\begin{proof} For simplicity of notations, we  drop the subscript $\lambda$ and divide the proof into five steps.
\medskip\\
(i) We first claim that it suffices to prove \eqref{Eq3} and \eqref{Eq39} for $p\leq q$. Indeed, suppose that $q<p$ and let
$$
r:=\tfrac{p}{p-1}<\theta:=\tfrac{q}{q-1}.
$$
By duality \eqref{Dual} and H\"older's inequality, we have
\begin{align*}
\|\nabla^2 u\|_{\mL^p_q(T)}&\stackrel{\eqref{Eq6}}{=}\sup_{g\in L^\infty_T(C^\infty_c), \|g\|_{\mL^r_\theta(T)}\leq 1}
\int^T_0\!\!\!\int_{\mR^d}\left(\int^t_0\cT_{s,t}f(s,x)\dif s\right) \nabla^2 g(t,x)\dif x\dif t\\
&=\sup_{g\in L^\infty_T(C^\infty_c), \|g\|_{\mL^r_\theta(T)}\leq 1}\int^T_0\!\!\!\int^t_0\left(\int_{\mR^d}\cT_{s,t}f(s,x)\nabla^2 g(t,x)\dif x\right)\dif s \dif t\\
&\stackrel{\eqref{Dual}}{=}\sup_{g\in L^\infty_T(C^\infty_c), \|g\|_{\mL^r_\theta(T)}\leq 1}\int^T_0\!\!\!\int^t_0\left(\int_{\mR^d}f(s,x)\cT^*_{s,t}\nabla^2 g(t,x)\dif x\right)\dif s \dif t\\
&=\sup_{g\in L^\infty_T(C^\infty_c), \|g\|_{\mL^r_\theta(T)}\leq 1}\int^T_0\!\!\!\int_{\mR^d}f(s,x)\left(\int^T_s\cT^*_{s,t}\nabla^2 g(t,x)\dif t \right)\dif x\dif s\\
&\leq C\sup_{g\in L^\infty_T(C^\infty_c), \|g\|_{\mL^r_\theta(T)}\leq 1}\|f\|_{\mL^p_q(T)}\|\nabla^2 g\|_{\mH^{-2,r}_\theta(T)}\leq C\|f\|_{\mL^p_q(T)},
\end{align*}
where the first inequality is due to \eqref{Eq39} for $p=r<\theta=q$.
\medskip\\
(ii) We only prove \eqref{Eq39} and \eqref{Eq301} for $p\leq q$ since \eqref{Eq3} and \eqref{Eq30} are similar.
By Marcinkiewicz's interpolation theorem (see \cite{St}), it suffices to prove that for any $p>1$ and $n\in\mN$,
\begin{align}\label{Eq4}
\left\|\nabla^2w\right\|_{\mH^{-2,p}_{np}(T)}\leq C\|f\|_{\mH^{-2,p}_{np}(T)}.
\end{align}
Below we fix $p>1$ and $n\in\mN$, and use the freezing coefficient argument to prove \eqref{Eq4}.
Let $\zeta$ be a nonnegative smooth function with support in the ball $B_\delta$ and
$\int_{\mR^d}\zeta^p\dif x=1$, where $\delta>0$ is a small constant and will be determined below. For $z\in\mR^d$, define
$$
\zeta_z(x):=\zeta(x-z),\ \ a_z(s):=a(s,z)
$$
and
$$
w_z(s,x):=w(s,x)\zeta_z(x),\ \ f_z(s,x):=f(s,x)\zeta_z(x).
$$
It is easy to see that
\begin{align}\label{Fr1}
\p_sw_z+\p_{ij} (a^{ij}_zw_z)-\lambda w_z+g_z=0,\ \ w_z(T)=0,
\end{align}
where
\begin{align*}
g_z:=f_z+\p_{ij} (a^{ij}w)\zeta_z-\p_{ij} (a^{ij}_zw\zeta_z).
\end{align*}
Moreover, by Fubini's theorem and $\int_{\mR^d}\zeta^p=1$, we have
\begin{align}\label{EG91}
\int_{\mR^d}\|w_z(s)\|_p^p\dif z=\int_{\mR^d}\|w(s)\zeta_z\|_p^p\dif z=\|w(s)\|_p^p.
\end{align}
Below we drop the time variable for simplicity. Noticing that
\begin{align*}
g_z=f\zeta_z-2\p_j(a^{ij}w)\p_i\zeta_z-a^{ij}w\p_{ij} \zeta_z+\p_{ij} ((a^{ij}-a^{ij}_z)w\zeta_z),
\end{align*}
and by Lemma \ref{Le21} with $\phi_z=\zeta_z,\p_i\zeta_z,\p_{ij}\zeta_z$ respectively,  we have
\begin{align}\label{HW8}
\begin{split}
\left(\int_{\mR^d}\|g_z\|^p_{-2,p}\dif z\right)^{1/p}
&\leq C\|f\|_{-2,p}+C_\delta\sum_{i,j}\|\p_j(a^{ij}w)\|_{-2,p}\\
&+C_\delta\sum_{i,j}\|a^{ij}w\|_{-2,p}+\omega_a(\delta)\|w\|_p,
\end{split}
\end{align}
where
$$
\omega_a(\delta):=\sup_{t\geq 0}\sup_{|x-y|\leq\delta}|a(t,x)-a(t,y)|.
$$
Let $a_n(t,x):=a(t,\cdot)*\rho_n(x)$ be the mollifying approximation of $a$. For every $\eps>0$, we can take $n$ large enough such that
\begin{align*}
&\sum_{i,j}\|\p_j(a^{ij}w)\|_{-2,p}+\sum_{i,j}\|a^{ij}w\|_{-2,p}\\
&\qquad\lesssim \|aw\|_{-1,p}\leq \|a_nw\|_{-1,p}+\|(a_n-a)w\|_{-1,p}\\
&\qquad\lesssim   \|a_n\|_{2,\infty} \|w\|_{-1,p}+\|(a_n-a)w\|_{p}\\
&\qquad\leq C_n\|w\|_{-1,p}+\omega_a(\tfrac{1}{n})\|w\|_{p}\\
&\qquad\lesssim \|w\|_{-2,p}+\eps\|w\|_{p},
\end{align*}
where the last step is due to the interpolation and Young's inequalities.
Hence, by \eqref{HW8}, for any $\eps\in(0,1)$ and $\delta>0$ being small enough,
\begin{align}\label{Eq7}
\left(\int_{\mR^d}\|g_z\|^p_{-2,p}\dif z\right)^{1/p}\lesssim \|f\|_{-2,p}+\|w\|_{-2,p}+\eps\|w\|_{p}.
\end{align}

\vspace{1mm}
\noindent(iii)
For any $s\in[0,T]$, notice that by Lemma \ref{Le21} again,
\begin{align}
\|\nabla^2 w\|^{np}_{\mH^{-2,p}_{np}(s,T)}&\lesssim\int^T_s\left(\int_{\mR^d} \|\nabla^2w(t)\zeta_z\|^p_{-2,p}\dif z\right)^n\dif t\no\\
&\lesssim\int^T_s\left(\int_{\mR^d} \|\nabla^2(w(t)\zeta_z)\|^p_{-2,p}\dif z\right)^n\dif t\no\\
&\quad+\int^T_s\left(\int_{\mR^d} \|\nabla w(t)\cdot\nabla\zeta_z\|^p_{-2,p}\dif z\right)^n\dif t\no\\
&\quad+\int^T_s\left(\int_{\mR^d} \|w(t)\cdot\nabla^2\zeta_z\|^p_{-2,p}\dif z\right)^n\dif t\no\\
&\lesssim\int^T_s\left(\int_{\mR^d} \|\nabla^2w_z(t)\|^p_{-2,p}\dif z\right)^n\dif t\no\\
&\quad+\int^T_s\|\nabla w(t)\|^{np}_{-2,p}\dif t+\int^T_s\|w(t)\|^{np}_{-2,p}\dif t\no\\
&\lesssim\int^T_s\!\!\!\!\int_{\mR^{nd}}\prod_{k=1}^n \|\nabla^2w_{z_k}(t)\|^p_{-2,p}\dif z_1\cdots\dif z_n\dif t\no\\
&\quad+\int^T_s\|w(t)\|^{np}_{-1,p}\dif t.\label{Eq5}
\end{align}
Given $z_1,\cdots, z_n\in\mR^d$ and by Lemma \ref{Le01}, we have
\begin{align*}
\int^T_s\prod_{k=1}^n\|\nabla^2w_{z_k}(t)\|_{-2,p}^p\dif t
\leq N\sum_{k=1}^n\int^T_s\|g_{z_k}(t)\|_{-2,p}^p\prod_{\ell\not=k}\|\nabla^2 w_{z_\ell}(t)\|_{-2,p}^p\dif t,
\end{align*}
which together with \eqref{Eq5} and \eqref{Eq7} yields that for any $\eps\in(0,1)$,
\begin{align*}
\|\nabla^2 w\|^{np}_{\mH^{-2,p}_{np}(s,T)}&\lesssim\sum_{k=1}^n\int^T_s\!\!\!\int_{\mR^{nd}}
\|g_{z_k}(t)\|_{-2,p}^p\prod_{\ell\not=k}\|w_{z_\ell}(t)\|_p^p\dif z_1\cdots\dif z_n\dif t+\|w\|^{np}_{\mH^{-1,p}_{np}(s,T)}\\
&=n\int^T_s\!\!\!\left(\int_{\mR^{d}}\|g_{z}(t)\|_{-2,p}^p\dif z\right)\!\!\!\left(\int_{\mR^d}\|w_z(t)\|_p^p\dif z\right)^{n-1}\dif t+\|w\|^{np}_{\mH^{-1,p}_{np}(s,T)}\\
&\!\!\!\stackrel{\eqref{EG91}}{=}n\int^T_s\left(\int_{\mR^{d}}\|g_{z}(t)\|_{-2,p}^p\dif z\right)\|w(t)\|_p^{(n-1)p}\dif t+\|w\|^{np}_{\mH^{-1,p}_{np}(s,T)}\\
&\!\!\!\stackrel{\eqref{Eq7}}{\lesssim}\|f\|^{np}_{\mH^{-2,p}_{np}(s,T)}+\|w\|^{np}_{\mH^{-2,p}_{np}(s,T)}+\eps\|\nabla^2w\|^{np}_{\mH^{-2,p}_{np}(s,T)},
\end{align*}
where the last step is due to H\"older's inequality and interpolation's inequality.
Taking $\eps=1/2$, we get for any $s\in[0,T]$,
\begin{align}\label{Eq8}
\|\nabla^2w\|^{np}_{\mH^{-2,p}_{np}(s,T)}\lesssim\|f\|^{np}_{\mH^{-2,p}_{np}(s,T)}+\|w\|^{np}_{\mH^{-2,p}_{np}(s,T)}.
\end{align}

\vspace{1mm}
\noindent
(iv) Let $A^z_{s,t}:=\int^t_s a_z(r)\dif r$ and
$$
P^z_{s,t} f(x):=\frac{1}{(2\pi)^{d/2}\det(A^z_{s,t})^{1/2}}\int_{\mR^d}\e^{-\<(A^z_{s,t})^{-1}y,y\>/2}f(x-y)\dif y.
$$
Notice that the solution of equation \eqref{Fr1} is explicitly given by
$$
w_z(s,x)=\int^T_s\e^{\lambda(s-t)}P^z_{s,t}g_z(t,x)\dif t.
$$
By \eqref{Non} and a standard interpolation technique, one sees that for any $\alpha\in[0,2)$,
there is a constant $C=C(\alpha,d,p,c_0)>0$ such that for all $z\in\mR^d$,
$$
\|w_z(s)\|_{\alpha-2,p}\leq C\int^T_s\frac{\e^{\lambda(s-t)}}{(t-s)^{\alpha/2}}\|g_z(t)\|_{-2,p}\dif t.
$$
Thus, for any $\alpha\in[0,2)$, by \eqref{GJ1} and Minkowski's inequality  we have
\begin{align}
&\|w(s)\|_{\alpha-2,p}\lesssim\left(\int_{\mR^d}\|w_z(s)\|^p_{\alpha-2,p}\dif z\right)^{\frac{1}{p}}
\leq\int^T_s\frac{\e^{\lambda(s-t)}}{(t-s)^{\alpha/2}}\left(\int_{\mR^d}\|g_z(t)\|_{-2,p}^p\dif z\right)^{\frac{1}{p}}\dif t\no\\
&\qquad\quad\stackrel{\eqref{Eq7}}{\leq}\int^T_s\frac{\e^{\lambda(s-t)}}{(t-s)^{\alpha/2}}\Big(\|f(t)\|_{-2,p}+\|w(t)\|_{-2,p}+\|\nabla^2w(t)\|_{-2,p}\Big)\dif t.\label{JH1}
\end{align}
Now by \eqref{JH1} with $\alpha=0$ and \eqref{Eq8} with $n=1$, we have
\begin{align*}
\|w(s)\|_{-2,p}^p
&\lesssim\int^T_s\Big(\|f(t)\|^p_{-2,p}+\|w(t)\|^p_{-2,p}+\|\nabla^2w\|^p_{-2,p}\Big)\dif t\\
&\lesssim\int^T_s\Big(\|f(t)\|^p_{-2,p}+\|w(t)\|^p_{-2,p}\Big)\dif t.
\end{align*}
which by Gronwall's inequality yields
$$
\|w\|^p_{\mH^{-2,p}_\infty(T)}=\sup_{s\in[0,T]}\|w(s)\|_{-2,p}^p\lesssim\|f\|^p_{\mH^{-2,p}_p(T)}\lesssim\|f\|^p_{\mH^{-2,p}_{np}(T)}.
$$
Substituting this into \eqref{Eq8} with $s=0$ and noting $\|w\|_{\mH^{-2,p}_{np}(T)}\lesssim \|w\|_{\mH^{-2,p}_\infty(T)}$, we obtain \eqref{Eq4}.
\medskip\\
(v) Finally, letting $q'=\frac{q}{q-1}$, for any $\alpha\in[0,2-\frac{2}{q})$, by \eqref{JH1} and H\"older's inequality, we have
\begin{align}
\|w(s)\|^q_{\alpha-2,p}&\lesssim\left(\int^T_s\frac{\e^{q'\lambda(s-t)}}{(t-s)^{\frac{q'\alpha}{2}}}\dif t\right)^{\frac{q}{q'}}\!\!\!\int^T_s\Big(\|f(t)\|_{-2,p}+\|w(t)\|_{-2,p}+\|w(t)\|_{p}\Big)^q\dif t\no\\
&\lesssim(1\vee\lambda)^{(\frac{\alpha}{2}-1+\frac{1}{q})q}\int^T_s\Big(\|f(t)\|^q_{-2,p}+\|w(t)\|^q_{-2,p}+\|\nabla^2w(t)\|_{-2,p}^q\Big)\dif t\no\\
&\stackrel{\eqref{Eq3}}{\lesssim}(1\vee\lambda)^{(\frac{\alpha}{2}-1+\frac{1}{q})q}\left(\|f\|^q_{\mH^{-2,p}_q(T)}+\int^T_s\|w(t)\|^q_{-2,p}\dif t\right),\label{HP9}
\end{align}
which yields by choosing $\alpha=0$ and Gronwall's inequality that
$$
\|w\|^q_{\mH^{-2,p}_\infty(T)}=\sup_{s\in[0,T]}\|w(s)\|_{-2,p}^q\lesssim\|f\|^q_{\mH^{-2,p}_q(T)}.
$$
The proof is complete by substituting this into \eqref{HP9}.
\end{proof}

\subsection{Proof of Theorem \ref{Th1}}
By standard continuity method (cf. \cite{Kry}), it suffices to establish the a priori estimate \eqref{Max}. We divide the proof into three steps.
\medskip\\
(i) (Case $b\equiv 0$) Fix $T>0$ and $p,q\in(1,\infty)$. Let $u\in \mH^{2,p}_q(T)$ and $f\in\mL^p_q(T)$ satisfy \eqref{PDE}.
Let $\rho_n$ be a family of mollifiers in $\mR^d$. Define
$$
u_n(t,x):=u(t,\cdot)*\rho_n(x),\ \ a_n(t,x):=a(t,\cdot)*\rho_n(x),\ \ f_n(t,x):=f(t,\cdot)*\rho_n(x).
$$
It is easy to see that $u_n$ satisfies
$$
\p_t u_n=a^{ij}_n\p_{ij} u_n-\lambda u_n+g_n,\ \ u_n(0)=0,
$$
where
$$
g_n:=f_n+(a^{ij}\p_{ij}u)*\rho_n-a^{ij}_n\p_{ij} u_n.
$$
Since $a_n$ satisfies {\bf (H$^a$)} uniformly in $n$ and $g_n\in\mL^q_T(H^{\infty,p})$,
for any $\alpha\in[0,2-\frac{2}{q})$, by \eqref{PDE11}, \eqref{Eq3} and \eqref{Eq30}, there is a $C>0$ such that for each $n\in\mN$ and $\lambda\geq 1$,
\begin{align*}
&\lambda^{1-\frac{\alpha}{2}-\frac{1}{q}}\|u_n\|_{\mH^{\alpha,p}_\infty(T)}+\|\p_t u_n\|_{\mL^p_q(T)}+\|\nabla^2 u_n\|_{\mL^p_q(T)}\\
&\quad\leq C\Big(\|f_n\|_{\mL^p_q(T)}+\|(a^{ij}\p_{ij}u)*\rho_n-a^{ij}_n\p_{ij} u_n\|_{\mL^p_q(T)}\Big).
\end{align*}
Letting $n\to\infty$ and by the property of convolutions, we obtain
\begin{align}\label{Eq31}
\lambda^{1-\frac{\alpha}{2}-\frac{1}{q}}\|u\|_{\mH^{\alpha,p}_\infty(T)}+\|\p_t u\|_{\mL^p_q(T)}+\|\nabla^2 u\|_{\mL^p_q(T)}\leq C\|f\|_{\mL^p_q(T)}.
\end{align}
Next, let $\chi^z_r$ be defined by \eqref{CHI}. Multiplying both sides of \eqref{PDE} by $\chi^z_r$, we have
$$
\p_t(u\chi^z_r)=a^{ij}\p_{ij}(u\chi^z_r)-\lambda u\chi^z_r+g^z_r,
$$
where
$$
g^z_r:=f\chi^z_r+\chi^z_r a^{ij}\p_{ij}u-a^{ij}\p_{ij}(u\chi^z_r).
$$
For any $\alpha\in[0,2-\frac{2}{q})$, by \eqref{Eq31} we have
\begin{align*}
\lambda^{1-\frac{\alpha}{2}-\frac{1}{q}}\|u\chi^z_r\|_{\mH^{\alpha,p}_\infty(T)}+\|\p_t u\chi^z_r\|_{\mL^p_q(T)}+\|\nabla^2 (u\chi^z_r)\|_{\mL^p_q(T)}\lesssim \|g^z_r\|_{\mL^p_q(T)}.
\end{align*}
Noticing that
$$
a^{ij}\p_{ij}(u\chi^z_r)-\chi^z_r a^{ij}\p_{ij}u=a^{ij}u\p_{ij}\chi^z_r+2a^{ij}\p_iu\p_j\chi^z_r,
$$
we have
$$
\|g^z_r\|_{\mL^p_q(T)}\lesssim \|f\chi^z_r\|_{\mL^p_q(T)}+\|u\chi^z_{2r}\|_{\mL^p_q(T)}
+\|\nabla u\cdot\chi^z_{2r}\|_{\mL^p_q(T)}.
$$
Hence, for any $\alpha\in[0,2-\frac{2}{q})$ and $\eps\in(0,1)$, by taking supremum in $z\in\mR^d$
and using \eqref{GW2}, we obtain that for all $\lambda\geq 1$,
\begin{align*}
&\lambda^{1-\frac{\alpha}{2}-\frac{1}{q}}\nor u\nor_{\widetilde\mH^{\alpha,p}_\infty(T)}+\nor \p_t u\nor_{\widetilde\mL^{p}_q(T)}+\nor u\nor_{\widetilde\mH^{2,p}_q(T)}\\
&\quad\lesssim\nor f\nor_{\widetilde\mL^p_q(T)}+\nor u\nor_{\widetilde\mL^p_q(T)}+\nor u\nor_{\widetilde\mH^{1,p}_q(T)}
\lesssim\nor f\nor_{\widetilde\mL^p_q(T)}+\nor u\nor_{\widetilde\mL^p_q(T)}+\eps\nor u\nor_{\widetilde\mH^{2,p}_q(T)},
\end{align*}
which implies by taking $\eps=1/2$ that
$$
\lambda^{1-\frac{\alpha}{2}-\frac{1}{q}}\nor u\nor_{\widetilde\mH^{\alpha,p}_\infty(T)}+\nor \p_t u\nor_{\widetilde\mL^{p}_q(T)}+\nor u\nor_{\widetilde\mH^{2,p}_q(T)}
\lesssim\nor f\nor_{\widetilde\mL^p_q(T)}+\nor u\nor_{\widetilde\mL^p_q(T)}.
$$
In particular, for $\alpha=0$, we have
$$
\nor u(T)\nor_p \lesssim\nor f\nor_{\widetilde\mL^p_q(T)}+\left(\int^T_0\nor u(s)\nor_{p}^q\dif s\right)^{1/q}.
$$
By Gronwall's inequality again, we obtain
$$
\nor u\nor_{\widetilde\mL^p_\infty(T)}\leq C\nor f\nor_{\widetilde\mL^p_q(T)},
$$
and so, for any $\alpha\in[0,2-\frac{2}{q})$,
\begin{align}\label{GK1}
\lambda^{1-\frac{\alpha}{2}-\frac{1}{q}}\nor u\nor_{\widetilde\mH^{\alpha,p}_\infty(T)}+\nor \p_t u\nor_{\widetilde\mL^{p}_q(T)}+\nor u\nor_{\widetilde\mH^{2,p}_q(T)}\lesssim\nor f\nor_{\widetilde\mL^p_q(T)}.
\end{align}
(ii) ($b\neq0$: subcritical case)
Let $q_1\in(\frac{2p}{p-d},q]$ and $\lambda\geq 1$.
For any $\alpha\in[0,2-\frac{2}{q_1})$, by \eqref{GK1}, we have
\begin{align}\label{GK2}
&\lambda^{1-\frac{\alpha}{2}-\frac{1}{q_1}}\nor u\nor_{\widetilde\mH^{\alpha,p}_\infty(T)}+\nor \p_t u\nor_{\widetilde\mL^{p}_q(T)}+\nor u\nor_{\widetilde\mH^{2,p}_{q_1}(T)}\no\\
&\quad\lesssim\nor f+b^i\p_i u\nor_{\widetilde\mL^p_{q_1}(T)}\leq \nor f\nor_{\widetilde\mL^p_{q_1}(T)}+\nor b^i\p_i u\nor_{\widetilde\mL^p_{q_1}(T)}.
\end{align}
Let $\frac{1}{q_2}+\frac{1}{q}=\frac{1}{q_1}$. For any $\theta\in(\frac{d}{p},1-\frac{2}{q_1})$, by H\"older's inequality and Sobolev's embedding \eqref{Sob}, we have
\begin{align}\label{GK4}
\nor b^i\p_i u\nor_{\widetilde\mL^p_{q_1}(T)}\leq \nor b\nor_{\widetilde\mL^{p}_{q}(T)}\nor u\nor_{\widetilde\mH^{1,\infty}_{q_2}(T)}
\lesssim \nor u\nor_{\widetilde\mH^{1+\theta,p}_{q_2}(T)}.
\end{align}
Substituting this into \eqref{GK2} with $\alpha=1+\theta$, we get
$$
\lambda^{\frac{1}{2}-\frac{\theta}{2}-\frac{1}{q_1}}\nor u\nor_{\widetilde\mH^{1+\theta,p}_\infty(T)}
\leq C\nor f\nor_{\widetilde\mL^p_{q_1}(T)}+ \nor u\nor_{\widetilde\mH^{1+\theta,p}_{q_2}(T)}.
$$
In particular, if $q_1<q$, then $q_2<\infty$ and by Gronwall's inequality again, we obtain
\begin{align}\label{GK5}
\nor u\nor_{\widetilde\mH^{1+\theta,p}_\infty(T)}\leq C\nor f\nor_{\widetilde\mL^p_{q_1}(T)}\leq C\nor f\nor_{\widetilde\mL^p_{q}(T)}.
\end{align}
The desired estimate now follows by \eqref{GK2}, \eqref{GK4} with $q_1=q$ and \eqref{GK5}.
\medskip\\
(iii) ($b\neq0$: critical case) Let $b_n(t,x):=b(t,\cdot)*\rho_{1/n}(x)$. Since $b\in \widetilde\mL^{d; {\rm uni}}_\infty$, by definition \eqref{GR1} we have
$$
\lim_{n\to\infty}\sup_{t\in[0,T]}\nor b_n(t)-b(t)\nor_d=0.
$$
Let $p<d$ and $q\in(1,\infty)$. For any $\eps\in(0,1)$, by Sobolev's embedding \eqref{Sob} and letting $n$ be large enough so that
$\sup_{t\in[0,T]}\nor b_n(t)-b(t)\nor_d\leq\eps$,
we have
\begin{align*}
\nor b^i\p_i u\nor_{\widetilde\mL^p_q(T)}&\leq \nor (b^i_n-b^i)\p_i u\nor_{\widetilde\mL^p_q(T)}+\nor b^i_n\p_i u\nor_{\widetilde\mL^p_q(T)}\\
&\leq \sup_{t\in[0,T]}\nor b_n(t)-b(t)\nor_d\nor\nabla u\nor_{\widetilde\mL^{pd/(d-p)}_{q}(T)}+\|b_n\|_\infty\nor u\nor_{\widetilde\mH^{1,p}_q(T)}\\
&\leq \eps\nor u\nor_{\widetilde\mH^{2,p}_{q}(T)}+C\|b_n\|_\infty\nor u\nor^{1/2}_{\widetilde\mL^{p}_q(T)}\nor u\nor^{1/2}_{\widetilde\mH^{2,p}_q(T)}\\
&\leq 2\eps\nor u\nor_{\widetilde\mH^{2,p}_{q}(T)}+C\|b_n\|_\infty^2\nor u\nor_{\widetilde\mL^{p}_q(T)}.
\end{align*}
Hence, for any $\alpha\in[0,2-\frac{2}{q})$, by \eqref{GK2} with $q_1=q$, we have
$$
\lambda^{1-\frac{\alpha}{2}-\frac{1}{q}}\nor u\nor_{\widetilde\mH^{\alpha,p}_\infty(T)}+\nor \p_t u\nor_{\widetilde\mL^{p}_q(T)}
+\nor u\nor_{\widetilde\mH^{2,p}_q(T)}\lesssim\nor f\nor_{\widetilde\mL^p_q(T)}
+\eps\nor u\nor_{\widetilde\mH^{2,p}_{q}(T)}+\nor u\nor_{\widetilde\mL^{p}_q(T)},
$$
which implies by taking $\eps=1/2$,
$$
\lambda^{1-\frac{\alpha}{2}-\frac{1}{q}}\nor u\nor_{\widetilde\mH^{\alpha,p}_\infty(T)}+\nor u\nor_{\widetilde\mH^{2,p}_q(T)}\lesssim\nor f\nor_{\widetilde\mL^p_q(T)}+
\nor u\nor_{\widetilde\mL^{p}_q(T)}.
$$
As above, by Gronwall's inequality, we obtain the desired estimate.

\section{Subcritical case: Proof of \autoref{Main1}}

In this section we assume {\bf (H$^\sigma$)} holds and for some $p_i,q_i\in[2,\infty)$ with $\frac{d}{p_i}+\frac{2}{q_i}<1$, $i=1,2$,
$$
\nabla\sigma\in \widetilde\mL^{p_1}_{q_1},\quad b\in \widetilde\mL^{p_2}_{q_2}.
$$
It is easy to see that {\bf (H$^a$)} holds for
$$
a^{ij}:=\sigma^{ik}\sigma^{jk}/2.
$$
We prepare the following crucial lemma for latter use.
\bl\label{Le41}
Let $X_t(x)$ be a solution of SDE \eqref{SDE1} and $p,q\in(1,\infty)$ with $\frac{d}{p}+\frac{2}{q}<2$.
\begin{enumerate}[(i)]
\item(Krylov's estimate) For any $T>0$, there is a constant $C>0$ such that for any $f\in \widetilde\mL^p_q(T)$ and $x\in\mR^d$,  $0\leq t_0<t_1\leq T$,
\begin{align}\label{PU2}
\bE\left(\int^{t_1}_{t_0} f(s,X_s(x))\dif s\Big|\sF_{t_0}\right)\leq C\nor f\nor_{\widetilde\mL^p_q(t_0,t_1)}.
\end{align}
\item (Khasminskii's estimate) For any $
\gamma\in\mR$ and $f\in  \widetilde\mL^p_q(T)$, we have
\begin{align}\label{Kas}
\bE\exp\left(
\gamma\int^T_0|f(s,X_s)|\dif s\right)<\infty.
\end{align}
\item (Generalized It\^o's formula) Let $p',q'\in[2,\infty)$ with $\frac{d}{p'}+\frac{2}{q'}<1$.
For any $u\in\widetilde\mH^{2,p'}_{q'}(T)$ with
$\p_tu\in\widetilde\mL^{p'}_{q'}(T)$, we have
\begin{align}\label{PU3}
\begin{split}
u(t,X_t)&=u(0,x)+\int^t_0(\p_s u+a^{ij}\p_i\p_j u+b^i\p_i u)(s,X_s)\dif s\\
&\qquad+\int^t_0(\sigma^{ij}\p_i u)(s,X_s)\dif W^j_s.
\end{split}
\end{align}
\end{enumerate}
\el
\begin{proof}
(i) By \eqref{Max} and using completely the same argument as in \cite[Theorem 5.7]{Xi-Zh},
we can prove the Krylov estimate \eqref{PU2}.
\medskip\\
(ii) Since $\frac{d}{p}+\frac{2}{q}<2$, we can choose $q'<q$ so that $\frac{d}{p}+\frac{2}{q'}<2$. Thus by \eqref{PU2} and H\"older's inequality we have
$$
\bE\left(\int^{t_1}_{t_0} f(s,X_s(x))\dif s\Big|\sF_{t_0}\right)\leq C\nor f\nor_{\widetilde\mL^p_{q'}(t_0,t_1)}\leq C(t_1-t_0)^{1-\frac{q'}{q}}\nor f\nor_{\widetilde\mL^p_{q}(T)},
$$
which implies \eqref{Kas} by \cite[Lemma 3.5]{Xi-Zh}.
\medskip\\
(iii) Let $u_n=(u*\rho_n)(t,x)$ be the mollifying approximation. By It\^o's formula we have
\begin{align}\label{II}
\begin{split}
u_n(t,X_t)&=u_n(0,X_0)+\int^t_0(\p_s u_n+a^{ij}\p_{ij} u_n+b^i\p_i u_n)(s,X_s)\dif s\\
&\quad+\int^t_0(\sigma^{ij}\p_i u_n)(s,X_s)\dif W^j_s.
\end{split}
\end{align}
For $R>0$, define a stopping time
$$
\tau_R:=\inf\{t\geq 0: |X_t|\geq R\}.
$$
Let $\chi_R$ be defined by \eqref{CHI}.
By It\^o's isometric formula, we have
\begin{align*}
&\bE\left|\int^{t\wedge\tau_R}_0(\sigma^{ij}\p_i(u_n-u))(s,X_s)\dif W^j_s\right|^2\\
&\qquad\leq\|\sigma\|^2_\infty\bE\left(\int^{t\wedge\tau_R}_0|\nabla (u_n-u)|^2(s,X_s)\dif s\right)\\
&\qquad \lesssim\bE\left(\int^{t}_0\chi^2_R(X_s)\cdot|\nabla (u_n-u)|^2(s,X_s)\dif s\right)\\
&\qquad \stackrel{\eqref{PU2}}{\lesssim}\nor \chi^2_R|\nabla (u_n-u)|^2\nor_{\widetilde\mL^{p'/2}_{q'/2}(T)}=\nor \chi_R\nabla (u_n-u)\nor^2_{\widetilde\mL^{p'}_{q'}(T)},
\end{align*}
which converges to zero by \eqref{LQ1} as $n\to\infty$. Similarly,
let $\frac{1}{p}:=\frac{1}{p_2}+\frac{1}{p'}$, $\frac{1}{q}:=\frac{1}{q_2}+\frac{1}{q'}$.
Since $\frac{d}{p}+\frac{2}{q}<2$, by \eqref{PU2} and H\"older's inequality we have
\begin{align*}
&\bE\left(\int^{t\wedge\tau_R}_0|b^i\p_i(u_n-u)|(s,X_s)\dif s\right)
\leq\bE\left(\int^t_0\chi_R(X_s)\cdot |b^i\p_i(u_n-u)|(s,X_s)\dif s\right)\\
&\quad\lesssim\nor \chi_R b^i\p_i(u_n-u)\nor_{\widetilde\mL^{p}_{q}(T)}
\leq \nor b\nor_{\widetilde\mL^{p_2}_{q_2}(T)}\nor\chi_{2R}\nabla(u_n-u)\nor_{\widetilde\mL^{p'}_{q'}(T)}\stackrel{n\to\infty}{\to} 0,
\end{align*}
and
\begin{align*}
\lim_{n\to\infty}\bE\left(\int^{t\wedge\tau_R}_0| (\p_s+a^{ij}\p_{i}\p_j)(u_n-u)|(s,X_s)\dif s\right)=0.
\end{align*}
By taking limits $n\to\infty$ for both sides of \eqref{II}, we get on $\{t\leq\tau_R\}$,
$$
u(t,X_t)=u(0,X_0)+\int^t_0(\p_s u+a^{ij}\p_i\p_j u+b^i\p_i u)(s,X_s)\dif s+\int^t_0(\sigma^{ij}\p_i u)(s,X_s)\dif W^j_s.
$$
Finally, letting $R\to\infty$, we obtain the desired formula.
\end{proof}

Below, we fix a $T>0$. Consider the following backward PDE:
$$
\p_t u+a^{ij}\p_{i}\p_j u-\lambda u+b^i\p_i u+b=0,\ \ u(T)=0.
$$
By \autoref{Th1}, there is a unique solution $u\in\widetilde\mH^{2,p_2}_{q_2}(T)$ such that for any $\alpha\in[0,2-\frac{2}{q_2})$ and $\lambda\geq 1$,
$$
\lambda^{1-\frac{\alpha}{2}-\frac{1}{q_2}}\nor u\nor_{\widetilde\mH^{\alpha,p_2}_\infty(T)}+\nor \p_tu\nor_{\widetilde\mL^{p_2}_{q_2}(T)}
+\nor u\nor_{\widetilde\mH^{2,p_2}_{q_2}(T)}\leq C\nor b\nor_{\widetilde\mL^{p_2}_{q_2}(T)}.
$$
In particular, since $\frac{d}{p_2}+\frac{2}{q_2}<1$, by \eqref{Sob} one can choose $\lambda$ large enough so that
\begin{align}\label{PU1}
\|u\|_\infty+\|\nabla u\|_\infty\leq \tfrac{1}{2}.
\end{align}
Define
$$
\Phi(t,x):=x+u(t,x).
$$
By \eqref{PU1}, one sees that $x\mapsto \Phi(t,x)$ is a $C^1$-diffeomorphism and
$$
\|\nabla\Phi\|_\infty,\ \ \|\nabla\Phi^{-1}\|_\infty\leq 2.
$$
Moreover, we also have
$$
\p_t \Phi+a^{ij}\p_{i}\p_j \Phi+b^i\p_i\Phi=\lambda u.
$$
Define
$$
\widetilde\sigma(t,y):=(\sigma^{ij}\p_i \Phi)(t,\Phi^{-1}(t,y))
$$
and
$$
\widetilde b(t,y):=\lambda u(t,\Phi^{-1}(t,y)).
$$
By the  generalized It\^o formula \eqref{PU3},  we have the following Zvonkin's transformation (see \cite[Theorem 3.10]{Xi-Zh}).

\bl\label{Le44}
$X_t$ solves SDE \eqref{SDE1} if and only if $Y_t=\Phi(t,X_t)$ solves the following SDE:
\begin{align}\label{New}
Y_t=y+\int^t_0\widetilde b(s,Y_s)\dif s+\int^t_0\widetilde \sigma(s,Y_s)\dif W_s\quad\text{with}\quad y:=\Phi(0,x).
\end{align}
\el

Now we can use the above lemma to prove \autoref{Main1}.

\begin{proof}[Proof of \autoref{Main1}]
By Lemma \ref{Le44}, it suffices to show the conclusions for SDE \eqref{New}.
Since the coefficients of SDE \eqref{New} are bounded and continuous, the existence of a solution $Y_t$ is well known.
By Yamada-Watanabe's theorem, we only need to prove the pathwise uniqueness for (\ref{New}) and show (i)-(iii)  for $Y$.

\vspace{1mm}
\noindent
(i) is proven in Lemma \ref{Le41}.

\vspace{1mm}
\noindent
(ii) For $i=1,2$, let $Y^{(i)}_t$ be two solutions of SDE \eqref{New} with starting point $y_i$, that is,
$$
Y^{(i)}_t=y_i+\int^t_0\widetilde b(s,Y^{(i)}_s)\dif s+\int^t_0\widetilde\sigma(s,Y^{(i)}_s)\dif W_s.
$$
For $\pp\geq 1$, by It\^o's formula we have
\begin{align}\label{LQ3}
|Y^{(1)}_t-Y^{(2)}_t|^{2\pp}=|y_1-y_2|^{2\pp}+\int^t_0|Y^{(1)}_s-Y^{(2)}_s|^{2\pp}\dif A_s+M_t,
\end{align}
where $M_t$ is a continuous local martingale given by
$$
M_t:=\int_0^t2\pp |Z_s|^{2\pp-2}\big[\widetilde\sigma(s,Y^{(1)}_s)-\widetilde\sigma(s,Y^{(2)}_s)\big]^*(Y^{(1)}_s-Y^{(2)}_s)\dif W_s,
$$
where the asterisk stands for the transpose of a matrix, and $A_t$ is defined by
\begin{align*}
A_t&:=\int^t_0\frac{2\pp\<Y^{(1)}_s-Y^{(2)}_s,\widetilde b(s,Y^{(1)}_s)-\widetilde b(s,Y^{(2)}_s)\>
+\pp\|\widetilde \sigma(s, Y^{(1)}_s)-\widetilde\sigma(s, Y^{(2)}_s)\|^2}{|Y^{(1)}_s-Y^{(2)}_s|^{2}}\dif s\\
&\quad+\int^t_0\frac{2\pp(\pp-1)|[\widetilde\sigma(s,Y^{(1)}_s)-\widetilde\sigma(s,Y^{(2)}_s)]^*(Y^{(1)}_s-Y^{(2)}_s)|^2}{|Y^{(1)}_s-Y^{(2)}_s|^{4}}\dif s.
\end{align*}
Notice that by Lemma \ref{Le2},
\begin{align*}
|\widetilde\sigma(s,x)-\widetilde\sigma(s,y)|
&\leq C|x-y|\Big(\cM_1|\nabla\widetilde\sigma(s,\cdot)|(x)+\cM_1|\nabla\widetilde\sigma(s,\cdot)|(y)+\|\widetilde\sigma\|_\infty\Big),\\
|\widetilde b(s,x)-\widetilde b(s,y)|
&\leq C|x-y|\Big(\cM_1|\nabla\widetilde b(s,\cdot)|(x)+\cM_1|\nabla\widetilde b(s,\cdot)|(y)+\|\widetilde b\|_\infty\Big).
\end{align*}
Thus, by the definitions of $\widetilde b$ and $\widetilde\sigma$ we have
\begin{align*}
|A_t|&\lesssim \int^t_0\left(\cM_1|\nabla\widetilde b|(s,Y^{(1)}_s)+\cM_1|\nabla\widetilde b|(s,Y^{(2)}_s)+\|\widetilde b\|_\infty\right)\dif s\\
&\quad+\int^t_0\left(\cM_1|\nabla\widetilde \sigma|^2(s,Y^{(1)}_s)+\cM_1|\nabla\widetilde \sigma|^2(s,Y^{(2)}_s)+\|\widetilde \sigma\|^2_\infty\right)\dif s\\
&\quad+\int^t_0\left(\cM_1|\nabla\widetilde \sigma|(s,Y^{(1)}_s)+\cM_1|\nabla\widetilde \sigma|(s,Y^{(2)}_s)+\|\widetilde \sigma\|_\infty\right)\dif s\\
&\lesssim t\left(\|\nabla\widetilde b\|_\infty+\|\widetilde b\|_\infty+\|\widetilde \sigma\|^2_\infty+\|\widetilde \sigma\|_\infty+1\right)\\
&\quad+\int^t_0\left(\cM_1|\nabla \sigma|^2(s,Y^{(1)}_s)+\cM_1|\nabla \sigma|^2(s,Y^{(2)}_s)\right)\dif s\\
&\quad+\int^t_0\left(\cM_1|\nabla^2 u|^2(s,Y^{(1)}_s)+\cM_1|\nabla^2 u|^2(s,Y^{(2)}_s)\right)\dif s,
\end{align*}
where we have used that $|\nabla \widetilde\sigma|(s,x)\lesssim |\nabla\sigma|(s,x)+|\nabla^2 u|(s,x)$.

On the other hand, by \eqref{GW1} we have
$$
\nor\cM_1 |\nabla \sigma|^2\nor_{\widetilde\mL^{p_1/2}_{q_1/2}(T)}\leq C \nor\,|\nabla\sigma|^2\nor_{\widetilde\mL^{p_1/2}_{q_1/2}(T)}
=C\nor \nabla\sigma\nor^2_{\widetilde\mL^{p_1}_{q_1}(T)}<\infty,
$$
and
$$
\nor\cM_1 |\nabla^2 u|^2\nor_{\widetilde\mL^{p_2/2}_{q_2/2}(T)}\leq C \nor\,|\nabla^2 u|^2\nor_{\widetilde\mL^{p_2/2}_{q_2/2}(T)}
=C\nor \nabla^2 u\nor^2_{\widetilde\mL^{p_2}_{q_2}(T)}<\infty.
$$
Thus, by Khasminskii's estimate \eqref{Kas},
$$
\bE\e^{\gamma A_T}<\infty, \quad \forall \gamma\in \mR.
$$
Hence, by \eqref{LQ3} and stochastic Gronwall's inequality (cf. \cite{Sc} or \cite[Lemma 3.7]{Xi-Zh}),
\begin{align}\label{wyy}
\bE\left(\sup_{t\in[0,T]} |Y^{(1)}_t-Y^{(2)}_t|^{\pp}\right)\leq C|y_1-y_2|^{\pp},
\end{align}
which in turn implies by \cite[Theorem 1.1]{XZ} that
$$
\sup_{y\in\mR^d}\bE\left(\sup_{t\in[0,T]}|\nabla Y_t(y)|^\pp\right)<\infty.
$$
Thus, by Lemma \ref{Le44} we obtain (\ref{Gr}). Moreover, by (\ref{wyy}) we also have the pathwise uniqueness.

\vspace{2mm}

\noindent(iii) Let $\widetilde \sigma_n(t,y):=\widetilde \sigma(t,\cdot)*\rho_n(y)$ be the usual mollifying approximation.
Let $Y_t^{n}$ be the unique strong solution of the following approximation SDE:
$$
\dif Y^{n}_t=\widetilde b(t,Y^{n}_t)\dif t+\widetilde \sigma_n(t, Y^{n}_t)\dif W_t, \ \ Y^{n}_0=y.
$$
By the classical Bismut-Elworthy-Li's formula (for example, see \cite{W-X-Zh}),  we have for any $h\in\mR^d$ and every bounded continuous function $\varphi$,
\begin{align}\label{Bis}
\nabla_h\bE\varphi\big(Y^{n}_t(y)\big)=\frac{1}{t}\bE\Bigg[\varphi\big(Y^{n}_t(y)\big)
\int_0^t\big[\widetilde\sigma_n\big(s,Y^{n}_s(y)\big)\big]^{-1}\nabla_h Y^{n}_s(y)\dif W_s\Bigg],
\end{align}
where $\nabla_h Y_t^{n}(y):=\lim_{\eps\to 0}[Y_t^{n}(y+\eps h)-Y_t^{n}(y)]/\eps$.  
On the other hand,
by ({\bf H$^\sigma$}) and the property of convolutions, it is easy to see that
$$
\lim_{|x-y|\to 0}\sup_n\sup_t\|\widetilde\sigma_n(t,x)-\widetilde\sigma_n(t,y)\|_{HS}=0,
$$
and for $n_0$ large enough,
$$
(2c_0)^{-1}|\xi|^2\leq|\widetilde\sigma_n(t,x)\xi|^2\leq 2c_0|\xi|^2,\ \ \xi\in\mR^d.
$$
Hence, $Y^{n}_t$ satisfies the Krylov estimate \eqref{PU2} with the constant $C$ independent of $n$.
As a result of \cite[Theorem 3.9]{Xi-Zh}, we have
$$
\lim_{n\to\infty}\bE\left(\sup_{t\in[0,T]} |Y^{n}_t(y)-Y_t(y)|\right)=0.
$$
Moreover,  as in the proof of \cite[(5.22)]{Zh4}, we have
$$
\lim_{n\to\infty}\sup_{y\in\mR^d}\bE\left(\sup_{t\in[0,T]} |\nabla Y^{n}_t(y)-\nabla Y_t(y)|\right)=0.
$$
Now taking limits $n\to\infty$ for both sides of \eqref{Bis} yields that for every $\varphi\in C^1_b(\mR^d)$,
$$
\nabla_h\bE\varphi\big(Y_t(y)\big)=\frac{1}{t}\bE\Bigg[\varphi\big(Y_t(y)\big)
\int_0^t\big[\widetilde\sigma\big(s,Y_s(y)\big)\big]^{-1}\nabla_h Y_s(y)\dif W_s\Bigg].
$$
Finally, using $\varphi\circ\Phi^{-1}_t(y)$ in place of $\varphi$ in the above formula, we obtain \eqref{feller}.
\end{proof}

\section{Critical case: Proof of \autoref{Main2}}
In this section we assume  that {\bf (H$^\sigma$)} holds and $b\in \widetilde\mL^{d; {\rm uni}}_\infty$.
Let
$$
b_n(t,x):=b(t,\cdot)*\rho_n(x),\ \ \sigma_n(t,x):=\sigma(t,\cdot)*\rho_n(x).
$$
By \eqref{GR1} and \eqref{LQ1}, it is easy to see that 
\begin{align}\label{KH0}
\sup_n\kappa^{b_n}_T(\eps)\leq C\kappa^{b}_T(\eps).
\end{align}
Without loss of generality we assume $s=0$ and consider the following approximation SDE:
$$
\dif X^n_t=b_n(t, X^n_t)\dif t+\sigma_n(t, X^n_t)\dif W_t, \ \ X^n_0=x.
$$
We first prove the following crucial lemma about Krylov's estimate.
\bl
Let $p\in(1,d)$ and $q\in(1,\infty)$ with $\frac{d}{p}+\frac{2}{q}<2$. For any $T>0$, there are constants $\theta=\theta(p,q)>0$
and $C>0$ such that for any $f\in C^\infty_c(\mR^{d+1})$, stopping time $\tau\leq T/2$ and $\delta\in(0,T/2)$,
\begin{align}\label{Kry0}
\sup_n\sup_{x\in\mR^d}\bE\left(\int^{\tau+\delta}_{\tau} f(s,X_s^n(x))\dif s\Big|\sF_{t_0}\right)\leq C\delta^\theta\nor f\nor_{\widetilde\mL^p_q(T)}.
\end{align}
\el
\begin{proof}
By discretizing stopping time approximation (see \cite[Remark 1.2]{Zh-Zh2}), it suffices to prove that  for any $0\leq t_0<t_1\leq T$ and $f\in C^\infty_c(\mR^{d+1})$.
\begin{align}\label{Kry11}
\sup_n\sup_{x\in\mR^d}\bE\left(\int^{t_1}_{t_0} f(s,X_s^n(x))\dif s\Big|\sF_{t_0}\right)\leq C(t_1-t_0)^\theta\nor f\nor_{\widetilde\mL^p_q(T)}.
\end{align}
Let $u_n$ be the smooth solution of the following backward PDE:
$$
\p_t u_n+\tfrac{1}{2}\sigma^{ik}_n\sigma^{jk}_n \p_{i}\p_ju_n+b^i_n\p_i u_n+f=0,\ u_n(t_1,\cdot)=0.
$$
Then, by It\^o's formula we have
$$
u_n(t_1,X^n_{t_1})=u_n(t_0,X^n_{t_0})-\int^{t_1}_{t_0}f(s,X^n_s)\dif s+\int^{t_1}_{t_0}\sigma^{ij}_n\p_i u_n(s,X^n_s)\dif W^j_s.
$$
Taking conditional expectation with respect to $\sF_{t_0}$, we obtain
$$
\bE\left(\int^{t_1}_{t_0}f(s,X^n_s)\dif s\Big|\sF_{t_0}\right)=u_n(t_0,X^n_{t_0})\leq\|u_n(t_0)\|_\infty.
$$
Since $\frac{d}{p}+\frac{2}{q}<2$, we can choose $q'<q$ so that $\frac{d}{p}+\frac{2}{q'}<2$.
Thus by \eqref{KH0}, \eqref{Max}, \eqref{Sob} and H\"older's inequality, there is constant $C>0$ such that
$$
\bE\left(\int^{t_1}_{t_0}f(s,X^n_{s})\dif s\Big|\sF_{t_0}\right)\leq C\nor f\nor_{\widetilde\mL^p_{q'}(t_0,t_1)}\leq C(t_1-t_0)^{1-\frac{q'}{q}}\nor f\nor_{\widetilde\mL^p_q(T)},
$$
which in turn gives  \eqref{Kry11}. The proof is complete.
\end{proof}

By the above lemma, we can show the following tightness result for $X^n$.
\bl\label{Le33}
For each $x\in\mR^d$, let $\mP^n_x$ be the law of $X^n_\cdot(x)$ in $\mC$. Then $(\mP^n_x)_{n\in\mN}$ is tight.
\el
\begin{proof}
Let $T>0$ and $\tau\leq T$ be any bounded stopping time. Notice that for every $\delta>0$,
$$
X^n_{\tau+\delta}-X^n_{\tau}=\int^{\tau+\delta}_\tau b_n(s,X^n_s)\dif s+\int^{\tau+\delta}_\tau\sigma_n(s, X_s)\dif W_s.
$$
Let $p\in(1,d)$ and $q\in(1,\infty)$ with $\frac{d}{p}+\frac{2}{q}<2$.
By  \eqref{Kry0} and Burkh\"older's inequality, there exists a $\theta>0$ such that for any $\delta\in(0,T)$,
\begin{align*}
\bE|X^n_{\tau+\delta}-X^n_{\tau}|&\leq\bE\left(\int^{\tau+\delta}_\tau |b_n(s,X^n_s)|\dif s\right)+C\bE\left(\int^{\tau+\delta}_\tau|\sigma_n(s, X_s)|^2\dif s\right)^{1/2}\\
&\leq C\delta^\theta\nor b_n\nor_{\widetilde\mL^p_q(2T)}+C\delta^{1/2}\stackrel{\eqref{LQ1}}{\leq} C\delta^\theta\nor b\nor_{\widetilde\mL^d_\infty(2T)}+C\delta^{1/2},
\end{align*}
where $C>0$ is independent of $n$.
Thus by \cite[Lemma 2.7]{ZZ}, we obtain
$$
\sup_n\bE\left(\sup_{s\in[0,T]}|X^n_{s+\delta}-X^n_s|^{1/2}\right)\leq C\left(\delta^{\theta/2}\nor b\nor_{\widetilde\mL^d_\infty(2T)}^{1/2}+\delta^{1/4}\right).
$$
By Chebyshev's inequality, we derive that for any $\eps>0$,
$$
\lim_{\delta\to 0}\sup_n\bP\left(\sup_{s\in[0,T]}|X^n_{s+\delta}-X^n_s|>\eps\right)=0,
$$
which implies the tightness of $X^n_\cdot$ by \cite[Theorem 1.3.2]{St-Va}.
\end{proof}

Now we can give the proof of Theorem \ref{Main2}.

\begin{proof}[Proof of \autoref{Main2}]
Since $(\mP^n_x)_{n\in\mN}\subset \sP(\mC)$ is tight, let $\mP_x$ be any accumulation point of $(\mP^n_x)_{n\in\mN}$.
By Krylov's estimate \eqref{Kry0}, it is by now easy to show that $\mP_x$
is a martingale solution of SDE \eqref{SDE1}, see for example, \cite{ZZ}. Moreover, \eqref{KR1} holds. 
We shall only prove  the uniqueness of martingale solutions.
Let $\mP_x^{(i)}\in\sM^{\sigma,b}_{0,x}, i=1,2$ be any two martingale solutions of SDE \eqref{SDE1} so that for any $T>0$, there is a constant $C>0$
such that for all $x\in\mR^d$ and $0\leq t_0<t_1\leq T$,  $f\in\widetilde\mL^p_q(t_0,t_1)$,
\begin{align}\label{Kr9}
\mE^{\mP^{(i)}_x}\left(\int^{t_1}_{t_0}f(s,\omega_s)\dif s\Big|\cB_{t_0}\right)\leq C\nor f\nor_{\widetilde\mL^p_q(t_0,t_1)}.
\end{align}
Let $p\in(1,d)$ and $q\in(1,\infty)$ satisfy $\frac{d}{p}+\frac{2}{q}<2$.
For $T>0$ and $f\in C^\infty_c([0,T]\times\mR^d)$, by Theorem \ref{Th1}, there is a unique solution $u\in\widetilde\mH^{2,p}_q(T)$
to the following backward equation:
$$
\p_t u+\sL^{\sigma,b}_tu+f=0,\ \ u(T)=0.
$$
Let $u_n(t,x):=u(t,\cdot)*\rho_n(x)$ be the mollifying approximation of $u$. Then we have
$$
\p_t u_n+\sL^{\sigma,b}_tu_n+g_n=0,\ \ u_n(T)=0,
$$
where
$$
g_n=f_n+(\sL^{\sigma,b}_tu)*\rho_n-\sL^{\sigma,b}_t(u*\rho_n).
$$
For $R>0$, define
$$
\tau_R:=\inf\{t\geq 0: |\omega_t|\geq R\}.
$$
By It\^o's formula, we have
\begin{align}\label{JH}
\mE^{\mP^{(i)}_x} u_n({T\wedge\tau_R}, \omega_{T\wedge\tau_R})=u_n(0,x)-\mE^{\mP^{(i)}_x} \left(\int^{T\wedge\tau_R}_0g_n(s,\omega_s)\dif s\right),\ \ i=1,2.
\end{align}
Since
$$
\nor\sL^{\sigma,b} u\nor_{\widetilde\mL^p_q(T)}\leq \|\sigma\|_\infty\nor\nabla^2u\nor_{\widetilde\mL^p_q(T)}
+\nor b\nor_{\widetilde\mL^{d}_\infty(T)}\cdot\nor\nabla u\nor_{\widetilde\mL^{pd/(d-p)}_q(T)}\stackrel{\eqref{Sob}}{\lesssim} \nor u\nor_{\widetilde\mH^{2,p}_q(T)},
$$
by Krylov's estimate \eqref{Kr9} and \eqref{LQ1}, we have
\begin{align*}
&\lim_{n\to\infty}\mE^{\mP^{(i)}_x} \left(\int^{T\wedge\tau_R}_0\Big((\sL^{\sigma,b}_tu)*\rho_n-\sL^{\sigma,b}_t(u*\rho_n)\Big)(s,\omega_s)\dif s\right)\\
&\quad\leq C\lim_{n\to\infty}\nor\chi_R((\sL^{\sigma,b}u)*\rho_n-\sL^{\sigma,b}(u*\rho_n))\nor_{\widetilde\mL^q_p(T)}=0,
\end{align*}
where the cutoff function $\chi_R$ is defined by \eqref{CHI}.
Letting $n\to\infty$ for both sides of  \eqref{JH} and by the dominated convergence theorem, we obtain
$$
\mE^{\mP^{(i)}_x} u({T\wedge\tau_R}, \omega_{T\wedge\tau_R})=u(0,x)-\mE^{\mP^{(i)}_x} \left(\int^{T\wedge\tau_R}_0f(s,\omega_s)\dif s\right),\ \ i=1,2,
$$
which, by letting $R\to\infty$ and noting $u(T)=0$, yields
$$
u(0,x)=\mE^{\mP^{(i)}_x} \left(\int^{T}_0f(s,\omega_s)\dif s\right),\ \ i=1,2.
$$
This in particular implies the uniqueness of martingale solutions (see \cite{St-Va}).
\end{proof}


\begin{thebibliography}{66}
\bib{Be-Fl-Gu-Ma}{article}{
   author={Beck, L.},
   author={Flandoli, F.},
   author={Gubinelli, M},
   author={Maurelli, M.},
   title={Stochastic ODEs and stochastic linear PDEs with critical drift: regularity, duality and uniqueness},
   eprint={ 1401.1530},
}
	
\bib{BL}{book}{
   author={J. Bergn},
   author={J. Lofstrom },
   title={Interpolation spaces},
   series={Grundlehren der math},
   volume={223},
   publisher={Springer-Verlag},
   date={1976},
  }

\bib{F-G-P}{article}{
   author={Flandoli, F.},
   author={Gubinelli, M.},
   author={Priola, E.},
   title={Well-posedness of the transport equation by stochastic
   perturbation},
   journal={Invent. Math.},
   volume={180},
   date={2010},
   number={1},
   pages={1--53},
}

\bib{Ki-Se}{article}{
    author={Kinzebulatov, D.},
   author={Semenov, Y. A.},
   title={Brownian motion with general drift},
   eprint={1710.06729v1},
}

	\bib{Ki}{article}{
   author={Kim, K. H.},
   title={$L^q(L^p)$-theory of parabolic PDEs with variable coefficients},
   journal={Bull. Korean Math. Society},
   volume={45 },
   date={2008},
   pages={169--190},
}


\bib{Kr1}{book}{
   author={Krylov, N. V.},
   title={Controlled diffusion processes},
   series={Translated from the Russian by A.B. Aries. Applications of Mathematics},
   volume={14},
   publisher={Springer-Verlag, New York-Berlin},
   date={1980},
  }
\bib{Kry}{book}{
   author={Krylov, N. V.},
   title={Lectures on elliptic and parabolic equations in H\"{o}lder spaces},
   series={Graduate Studies in Mathematics},
   volume={12},
   publisher={American Mathematical Society, Providence, RI},
   date={1996},
  }

\bib{K3}{article}{
   author={Krylov, N. V.},
   title={The heat equation in $L_q((0,T), L_p)$-spaces with weights},
   journal={SIAM Journal Math. Anal},
   volume={32},
   date={2001},
   pages={1117--1141},
}

\bib{Kr-Ro}{article}{
   author={Krylov, N. V.},
   author={R\"{o}ckner, M.},
   title={Strong solutions of stochastic equations with singular time
   dependent drift},
   journal={Probab. Theory Related Fields},
   volume={131},
   date={2005},
   number={2},
   pages={154--196},
}
\bib{M-N-P-Z}{article}{
   author={Menoukeu, P. O.},
   author={Meyer, B. T.},
   author={Nilssen, T.},
   author={Proske, F.},
   author={Zhang, T.},
   title={A variational approach to the construction and Malliavin
   differentiability of strong solutions of SDE's},
   journal={Math. Ann.},
   volume={357},
   date={2013},
   number={2},
   pages={761--799},
}
\bib{Mo-Ni-Pr}{article}{
   author={Mohammed, S. E. A.},
   author={Nilssen, T.},
   author={Proske, F.},
   title={Sobolev differentiable stochastic flows of SDE's with singular coefficients},
   journal={Ann. Prob.},
   volume={43},
   date={2015},
   pages={1535--1576},
}
\bib{Na}{article}{
    author={Nam, K.},
   title={Stochastic differential equations with critical drifts},
   eprint={1802.00074},
}
\bib{Sc}{article}{
   author={Scheutzow, M.},
   title={A stochastic Gronwall's lemma},
   journal={Infinite Dimensional Analysis, Quantum Probability and Related Topics},
   volume={16},
   date={2013},
}
\bib{St}{book}{
   author={Stein, E. M.},
   title={Singular integrals and differentiability properties of functions},
   series={Princeton Mathematical Series },
   volume={30},
   publisher={Princeton University Press, Princeton, NJ},
   date={1970},
  }

  \bib{St-Va}{book}{
   author={Stroock, D. W.},
   author={ Varadhan, S. R. S.},
   title={ Multidimensional diffusion processes},
   series={Grundlehren der Mathematischen Wissenschaften},
   volume={233},
   publisher={Springer-Verlag, Berlin-New York},
   date={1979},
  }
  \bib{Ve}{article}{
   author={Veretennikov, A. J.},
   title={On the strong solutions of stochastic differential equations},
   journal={Theory Probab. Appl.},
   volume={24},
   date={1979},
   pages={354--366},
}

\bib{W-X-Zh}{article}{
   author={Wang, L.},
   author={ Xie, L. },
   author={ Zhang, X.},
   title={Derivative formulae for SDEs driven by multiplicative $\alpha$-stable-like processes},
   journal={Stoch. Proc. Appl.},
   volume={125},
   date={2015},
   pages={867--885},
}

\bib{XZ}{article}{
   author={Xie, L.},
   author={Zhang, X.},
   title={Sobolev differentiable flows of SDEs with local Sobolev and
   super-linear growth coefficients},
   journal={Ann. Probab.},
   volume={44},
   date={2016},
   pages={3661--3687},
}
		
\bib{Xi-Zh}{article}{
   author={Xie, L.},
   author={Zhang, X.},
   title={Ergodicity of stochastic differential equations with jumps and singular coefficients},
  journal={{\rm To appear in} Annales de l'Institut Henri Poincar\'e - Probabilit\'es et Statistiques},
   eprint={1705.01402v1},
}
\bib{Zh3}{article}{
   author={Zhang, X.},
   title={Strong solutions of SDES with singular drift and Sobolev diffusion
   coefficients},
   journal={Stochastic Process. Appl.},
   volume={115},
   date={2005},
   pages={1805--1818},
}
	
\bib{Zh1}{article}{
   author={Zhang, X.},
   title={Stochastic homemomorphism flows of SDEs with singular drifts and Sobolev diffusion coefficients},
   journal={Electron. J. Probab.},
   volume={16},
   date={2011},
   pages={1096--1116},
}

\bib{Zh4}{article}{
   author={Zhang, X.},
   title={Stochastic differential equations with Sobolev coefficients and applications},
   journal={Annals of Applied Prob.},
   volume={26},
   date={2016},
   pages={2697--2732},
}
\bib{ZZ}{article}{
   author={ Zhang, X.},
   author={Zhao, G.},
   title={Heat kernel and ergodicity of SDEs with distributional drifts},
   eprint={1710.10537v2},
}
\bib{Zh-Zh2}{article}{
   author={ Zhang, X.},
   author={Zhao, G.},
   title={Stochastic Lagrangian path for Leray solutions of 3D Navier-Stokes
  equations},
   eprint={1904.04387},
}




\end{thebibliography}
\end{document}